\documentclass[12pt]{article}
\begin{document}

\newtheorem{lemma}{Lemma}
\newtheorem{theorem}[lemma]{Theorem}
\newtheorem{Remark}[lemma]{Remark}
\newtheorem{Fact}[lemma]{Fact}
\newenvironment{proof}{P r o o f.}{\rule{2mm}{3mm}\vspace{5mm}}
\def\rng{{\rm rng}}
\def\dom{{\rm dom}}
\def\text{\rm}

\title {Balanced $d$-lattices are complemented\thanks{This paper
is also available from {\tt www.arxiv.org} and from authors' homepages}}
\date{}
\author{Martin Goldstern\\
{\it Algebra, TU Wien, Wiedner Hauptstrasse 8-10/118,}\\
{\it A-1040 Wien, Austria}\\
{\small
\tt martin.goldstern@tuwien.ac.at, http://www.tuwien.ac.at/goldstern/}\\
 Miroslav Plo\v s\v cica
\thanks{Supported by VEGA Grant 2/1131/21.}\\
{\it Mathematical Institute, Slovak Academy of Sciences,
Gre\v s\'akova 6,}\\
{\it 04001 Ko\v sice, Slovakia}\\
{\small
\tt ploscica@saske.sk, http://www.saske.sk/MI/eng/ploscica.htm}}

\maketitle


\vskip0.5cm
According to Chajda and Eigenthaler (\cite{ce}),
a $d$-lattice is a bounded lattice $L$ satisfying
for all $a,c\in L$ the implications
\begin{description}
\item[(i)] $(a,1)\in\theta(0,c)\ \to\ a\vee c=1$;
\item[(ii)] $(a,0)\in\theta(1,c)\ \to\ a\wedge c=0$;
\end{description}
where $\theta(x,y)$ denotes the least congruence on $L$ containing
the pair $(x,y)$. Every bounded distributive lattice is a $d$-lattice.
The 5-element nonmodular lattice $N_5$ is a $d$-lattice.

\begin{theorem} A bounded lattice is a $d$-lattice if and only if
all maximal ideals and maximal filters are prime.\label{charact}
\end{theorem}
\begin{proof}
Let $I$ be a maximal ideal in a $d$-lattice $L$. Let $x,y\in L\setminus I$.
We need to show that $x\wedge y\in L\setminus I$. Since $I$ is maximal,
there are $c_1, c_2\in I$ such that $c_1\vee x=c_2\vee y=1$. For
$c=c_1\vee c_2\in I$ we have $c\vee x=c\vee y=1$. Then $(x,1)=
(0\vee x,c\vee x)\in\theta(0,c)$ and similarly $(y,1)\in\theta(0,c)$,
hence $(x\wedge y, 1)\in\theta(0,c)$. By (i) we have $(x\wedge y)\vee c=1$,
hence $x\wedge y\not\in I$. The primality of maximal filters can be proved
similarly.

Conversely, assume that all maximal ideals and filters in $L$ are prime.
To show (i), assume that $a,c\in L$, $a\vee c\ne 1$. By the Zorn lemma,
there exists a maximal ideal $I$ containing $a\vee c$. By our assumption,
$I$ is prime. Then $\alpha=I^2\cup(L\setminus I)^2$ is a congruence
on $L$. Since $c\in I$, we have $(0,c)\in\alpha$, which implies that
$\theta(0,c)\subseteq\alpha$. Since $a\in I$, we have $(a,1)\not\in\alpha$,
hence $(a,1)\not\in\theta(0,c)$. This shows (i). The proof of (ii)
is similar.
\end{proof}

By \cite{ce}, a bounded lattice is called ``balanced'',
if the $0$-class of any
congruence determines the 1-class, and conversely.  They showed that
complemented  lattices are balanced, and they asked:
 \begin{itemize}
\item [$(*)$]
Is there a $d$-lattice which is balanced but  not
complemented?
\end{itemize}

We use the above characterization of $d$-lattices to answer this question.

If $A$ is a subset of an algebra,
write $\theta_A $ for the smallest congruence that identifies all
elements of~$A$; if $\phi$ is a congruence, $x$ an element, write
$x/\phi$ for the $\phi$-congruence class of~$x$.

Further,
a congruence $\phi$ (on an algebra with constants $0$ and
$1$) is called balanced if $ 0/\phi = 0/\theta_{(1/\phi)}$ and
 $ 1/\phi = 1/\theta_{(0/\phi)}$;
an algebra  is called balanced iff all
 its congruence relations are balanced, or equivalently if:  for
 any congruence relations~$\phi$, $\phi'$ we have:
 \begin{quote}
 $0/\phi = 0 / \phi' $ iff
 $1/\phi = 1 / \phi' $.
 \end{quote}

Fix a $d$-lattice $(L,\vee,\wedge,0,1)$.
For $a\in L$ we denote $F_a:= \{x: x\vee a = 1\}$, and $I_a:=
\{x: x\wedge a = 0 \}$.

\begin{Fact} $F_a$ is a filter, $I_a$ is an ideal.
\end{Fact}

\begin{proof} Let $x,y\in F_a$. Similarly as in the proof of
Theorem \ref{charact}, $(x,1)\in\theta(0,a)$, $(y,1)\in\theta(0,a)$,
hence $(x\wedge y,1)\in\theta(0,a)$, which by the definition of
a $d$-lattice implies $x\wedge y\in F_a$. The proof for $I_a$ is
similar.
\end{proof}

\begin{Fact}\label{dis}
If $I$ is an ideal disjoint to~$F_a$, and $a\notin I$, then
also the ideal generated by $I\cup \{a\}$ is disjoint to~$F_a$.
\end{Fact}

\begin{proof} If $x\le i \vee a$ for some $i\in I$, and $x\in F_a$, then
also $i\vee a\in F_a$, hence $i\vee a=(i\vee a)\vee a=1$. Thus, $i\in
F_a$, so $F_a\cap I\not=\emptyset$.
\end{proof}

\begin{Fact}\label{level}
If $f:L_1 \to L_2$ is a homomorphism from $L_1$ onto
$L_2$, and $L_1$ is balanced, then $L_2$ is balanced.
\end{Fact}

\begin{proof} In fact, this holds ``level-by-level'':  If $\phi$ is an
unbalanced congruence on~$L_2$, then the preimage of $\phi$ is
unbalanced on~$L_1$.
\end{proof}

\begin{theorem}
The following are equivalent (for a $d$-lattice~$L$):
\begin{enumerate}
\item There is a maximal (hence prime)  filter whose complement
is not a maximal ideal.
\item There is a maximal (hence prime)  ideal whose complement
is not a maximal filter.
\item There are two prime ideals in~$L$, one properly containing the
other.
\item There are two prime filters in~$L$, one properly containing the
other.
\item There is a homomorphism from $L$ onto the 3-element lattice $\{0,d,1\}$.
\item $L$ is not balanced.
\item $L$ is not complemented.
\end{enumerate}

In particular a $d$-lattice is balanced iff it is
complemented.
\end{theorem}

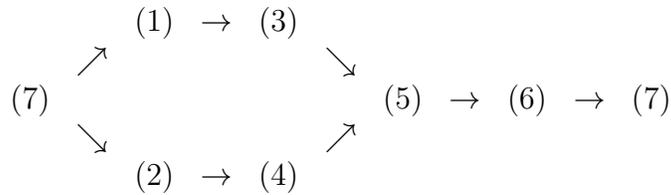
\begin{proof}
$$
\begin{array}{ccccccccccc}
    &        & \smash{(1)}& \to & (3) &\\
    &\nearrow& &&   &\searrow\\
(7) &        & &&   &            & (5) &  \to   & (6) &  \to   & (7) \\
    &\searrow& &&   &\nearrow&\\
    &        & (2)& \to & (4)&

\end{array}
$$
(1) $\to $ (3): By \ref{charact}, the complement of a maximal filter is
a (necessarily prime) ideal. If this ideal is not maximal, it can be properly
extended to a maximal (hence prime) ideal. The proof of (2) $\to$ (4)
is similar (dual).

(3) $\to $ (5):  Let $I_1\subset I_2\subset L$ be prime ideals.
Map $I_1$ to~$0$,
$I_2\setminus I_1$ to~$d$, and $L\setminus I_2$ to~$1$.  Check that
this is a lattice homomorphism.
The proof of (4) $\to $ (5) is dual.

(5) $\to $ (6) follows from fact \ref{level}, since the three-element
  lattice is not balanced.

(6) $\to $ (7) is from \cite{ce}.

\bigskip
Now we show (7) $\to $ (1).   (Again, (7) $\to $ (2) is dual.)
Assume that $L$ is not complemented, so there is some $a$
such that
 $F_a\cap I_a=\emptyset$.
Let $F_1$ be the filter generated by $F_a\cup \{a\}$.  We have
$F_1\cap I_a=\emptyset$ by the dual of Fact \ref{dis}, so $F_1$
is proper. By the Zorn lemma,
$F_1$ can be extended to a maximal filter $F$. Let $I_1=L\setminus F$.
It is enough to see that $I_1$ is not maximal.
Let $I$ be the ideal generated by $I_1\cup\{a\}$.
By Fact \ref{dis},  $I \cap F_a = \emptyset$, so $I$ is a proper ideal
properly extending $I_1$.
\end{proof}

\end{document}